\documentclass[12pt]{article}
\usepackage{graphics, latexsym}
\usepackage{graphicx}
\usepackage{setspace}
\usepackage{geometry}
\usepackage{amssymb}
\usepackage{amsmath}
\usepackage{float}
\usepackage{authblk}
\usepackage{booktabs}
\usepackage{threeparttable}
\usepackage{multirow}
\usepackage{diagbox}
\usepackage{color}
\newtheorem{definition}{Definition}
\newtheorem{theorem}{Theorem}

\newtheorem{remark}{Remark}

\begin{document}
\begin{spacing}{2.0}
\title{Optimal Sequential Tests for Detection of Changes under Finite Measure Space for Finite Sequences of Networks}
\author[1]{Lei Qiao\^*}
\author[1]{Dong Han}
\affil[1]{Department of Statistics, Shanghai Jiao Tong University,
  China.}

\maketitle

\begin{abstract}
This paper considers the change-point problem for finite sequences of networks. To avoid the difficulty of computing the normalization coefficient in the models such as Exponential Random Graphical Model (ERGM) and Markov networks, we construct a finite measure space with measure ratio statistics. A new performance measure of detection delay is proposed to detect the changes in distribution of the network data. And under the performance measure we defined, an optimal sequential test is presented. The good performance of the optimal sequential test is illustrated numerically on ERGM and Erd\H{o}s-R\'{e}nyi network sequences.
\end{abstract}

\renewcommand{\thefootnote}{\fnsymbol{footnote}}
\footnotetext{
Supported by RGC Competitive Earmarked Research Grants, National Basic Research Program of China
(973 Program, 2015CB856004) and the National Natural Science Foundation of China (11531001)
\newline $\,\, MSC\, 2010\,subject\,classification $. Primary 62L10;
Secondary 62L15
\newline *corresponding author
\newline Email address:ql0130@sjtu.edu.cn(Lei Qiao), donghan@sjtu.edu.cn(Dong Han)}
\emph{Keywords}: ERGM; optimal sequential test; finite measure space; change-point detection; finite sequences of networks.

\section{Introduction}
One of the basic problems of quickest change-point detection is the design of an optimal sequential test. If the alarm is too fast, there is a high risk that the detection is false. At the same time, if the alarm is too slow, the detection delay would be too large. So an optimal sequential test is important, which is expected to have the smallest detection delay among all tests without exceeding a certain false alarm rate. 

Optimal sequential test has vast applications in many fields, such as industrial quality and process control, information and communication systems, biostatistics, clinical trials and public-health, econometrics and financial surveillance, and so on. With rapid development in technology, network data are often encountered, such as social networks, earthquake networks, traffic networks and so on. Due to the increasing amount of network data, it becomes important and useful to analyze network data, especially the change detection problem for observation sequences of network data, which has been shown in Amini and Nguyen(2013); Keshavarz, Michailidis, and Atchade (2018); Wei, Wiesel, and Blum (2012). These papers discussed the change detection problem in Gaussian graphical model. For Gaussian graphical model, accurate distributions can be conveniently obtained, which makes the change-point problem easy to deal with. Thus, Gaussian graphical model is a very special case. 

However, there are many challenges and difficulties in dealing with large scale network data. When the network is large, it is difficult to obtain the distribution of the network data because the normalization coefficient is not easy to be calculated and needs exponential computational complexity. For example, the distribution of Exponential Random Graphical Model(ERGM), also known as $p^*$ model, introduced by Frank and Strauss(1986); Frank(1991); Wasserman and Pattison(1996), is 
$$
P(X=x)=\exp[\theta^T h(x)][Z(\theta)]^{-1},
$$
where $x$ represents a network, $h(x)$ are the features of the network, and the normalization coefficient denoted by $Z(\theta)$ is difficult to obtain especially for large scale network data. ERGM model is popular in recent years, see Vu, Hunter, and Schweinberger (2013); Hanneke, Fu, and Xing (2010); Krivitsky and Handcock (2014). These researches used the ERGM or extended ERGM models for statistical analysis of network data. Using ERGM is the practical requirement. And this phenomenon also exists for the Markov networks of probabilistic graphical model introduced by Koller and Friedman (2009). Hence, it urges us to find a method to overcome this difficulty. 

There are mainly four kinds of optimal sequential tests in the literature: Shiryaev test $T_S$ (Shiryaev (1963, 1978)), SLR (sum of the log likelihood ratio) test $T_{SLR}$ (Chow, Robbins and Siegmund (1971, 108); Fris\'{e}n (2003)), conventional CUSUM test $T_C$ (Page (1954); Moustakides (1986)) and Shiryaev-Roberts test $T^r_{SR}$ (Polunchenko and Tartakovsky (2010)). A common ground of the four kinds of optimal sequential tests is that they are based on likelihood ratio statistics, and the analyzed observation sequences are one-dimensional infinite and independent data. However, it is unrealistic to obtain infinite sequences of networks. Thus, there is a necessity to deal with the problem of how to give the corresponding optimal sequential test for finite sequences of networks. 

In this paper, in order to overcome the above difficulties, i.e., the calculation of normalization coefficient and the finite sequences of networks, we construct a finite measure space to replace the traditional probability space by discarding the normalization coefficient $Z(\theta)$. And under the finite measure space, we propose a new measure of detection delay. The likelihood ratio statistics are replaced by measure ratio statistics. In a word, the method we proposed is based on the finite measure space rather than the probability space. We know that the difference between the finite measure space and the probability space is the normalization coefficient, so the traditional likelihood ratio statistic is a special case of our method. Moreover, we consider unweighted network sequences here. A network is a set of vertices (nodes) with connections (edges) between them. If no weight is assigned to each edge, we only have connected or disconnected information, then the network is called unweighted network. Furthermore, the optimal sequential test in this paper is applied to dependent and finite network observation sequences. 

The rest of this paper is organized as follows. Section \ref{section2} states the method we use in this paper. Section \ref{section2.1} presents a performance measure, $\mathcal{J}_{N}(T)$, to evaluate the performance of a sequential test in detecting changes in finite network observation sequences. Section \ref{section2.2} gives the optimal sequential test under the performance measure we defined. In Section \ref{section3}, we provide several numerical simulations to show that the optimal sequential test we proposed has better performance than other three sequential tests. Section \ref{section3.1} and Section \ref{section3.2} are the cases of independent and dependent finite sequences of networks in the finite measure space, respectively. Section \ref{section3.3} and Section \ref{section3.4} are the cases of independent and dependent finite sequences of networks in the traditional probability space, respectively. Section \ref{section4} includes some concluding remarks. Proofs of the theorems are given in the Appendix.

\section{Optimal sequential tests for finite sequences of networks under the finite measure space} \label{section2}
In this section, we first briefly describe the difference between the probability space and the finite measure space, then construct the optimal sequential tests for $N$ network observations.

In this paper, we consider the models such as the Exponential Random Graphical model (ERGM) proposed in Frank and Strauss(1986); Frank(1991); Wasserman and Pattison(1996). The distribution is
$$
P(X=x)=\frac{\exp(\theta^T h(x))}{Z(\theta)},
$$
where $x$ represents a network, $h(x)$ are the crucial features of the network (e.g. degree, betweenness, closeness, eigenvector centrality, node eccentricity, node strength and node domination) and $Z(\theta)$ is the normalization coefficient which requires complicated calculations. There are many different features of networks, which are reflected in $h(x)$, a multi-dimensional vector.

Since the normalization coefficient is difficult to obtain, the statistics based on likelihood ratio in the probability space, such as CUSUM-type, cannot be used. To avoid this, we consider a finite measure space $(\Omega, \mathcal{F}, \textbf{M})$. We extend the likelihood ratio statistics to the measure ratio statistics, and replace the probability space with a finite measure space.

As with the probabiliy space, we also have measure density functions and conditional measure density functions for network objects. The definitions are shown below:
\begin{itemize}
	\item [(i)] Suppose that there is a network set $\Omega$, let $\textbf{M}$ be a finite measure, assume that $\textbf{M}$ is differentiable in the usual sense in calculus. Let $m$ be the derivative of $\textbf{M}$,
	\begin{equation}
	\label{section2:measure}
	\textbf{M}(A)=\int_A m(x)dx \,\,\,\,\,A \in \mathcal{F},
	\end{equation}
	$x$ is the network object, $\mathcal{F}$ is $\sigma$-field. In the discrete case, we use $\sum$ to replace $\int$. Then $m$ is the corresponding measure density function. 
	\item [(ii)] The conditional measure density function is defined as follows,
	$$
	m(x|y) = \frac{m(x,y)}{m(y)}.
	$$
\end{itemize}

We consider unweighted network sequences here. The definition of the joint measure density function is similar, which will be explained in detail later. Suppose that there is a finite and dependent network observation sequence, denoted by $\{X_k, 0 \leq k \leq N \}$. Let $d$ denote the number of nodes in each network. Then, we consider the change-point problem for network sequences as follows:

\begin{itemize}
	\item [(i)]All networks in the sequence have a common number of nodes $d$, which is a constant. Without loss of generality, we assume $N\geq 2$.
	
	\item [(ii)]The pre-change and post-change measure density functions are denoted by $m_0$ and $m_1$, and conditional measure density functions of the $j$-th observation $X_j$ for $1\leq j\leq N$ are denoted by $m_{0j}$ and $m_{1j}$, which are all known. The measures above are all finite. 	
\end{itemize}

Let $\tau=k$ ($1 \leq k \leq N $) be the change-point, the joint measure density function $m_k$ can be written as
\begin{equation}
m_k(X_0, X_1, ..., X_N)=m_{0}(X_0)\prod_{j=1}^{k-1}m_{0j}(X_j|X_{j-1}, ..., X_0)\prod_{j=k}^N m_{1j}(X_j|X_{j-1}, ..., X_0),
\end{equation}
where $m_{0j}(X_j|X_{j-1}, ..., X_0)$ and $m_{1j}(X_j|X_{j-1}, ..., X_0)$ are respectively the pre-change and the post-change conditional measure density functions of the $j$-th observation $X_j$ for $1\leq j\leq N$. When $\tau=k >N$, $\prod_{j=k}^N:=1$, i.e., a change never occurs in $N$ observations $X_1, X_2, ..., X_N$. The corresponding joint measure function is denoted by $\textbf{M}_k$, which will be introduced in detail later.

\subsection{Performance measures of sequential tests} \label{section2.1}
In general, before considering the optimality of a sequential test, we should introduce a performance measure to evaluate how well a sequential test detects changes in observation sequences. There are some different kinds of performance measures for sequential tests which were used in Shiryaev (1963, 1978); Chow, Robbins, and Siegmund (1971, 108); Fris\'{e}n (2003); Moustakides (1986); Polunchenko and Tartakovsky (2010); Han, Tsung, and Xian (2017). These performance measures rely on the probability distribution $\textbf{P}_k$ or the probability expectation $\textbf{E}_k$ in the probability space. Similarly, in order to propose a new performance measure in the finite measure space, we first need to give the definitions about the finite measure $\textbf{M}_k$ and the measure expectation $\textbf{ME}_k$ when the change-point $\tau=k$ ($1 \leq k \leq N $). Naturally, the new performance measure under the finite measure space relies on $\textbf{M}_k$ or $\textbf{ME}_k$ ($1 \leq k \leq N$).

We give the definitions about the finite measure $\textbf{M}_k$ and the measure expectation $\textbf{ME}_k$ in the following when the change-point $\tau=k$ ($1 \leq k \leq N $). Specifically, the definition about $\textbf{M}_k$ is shown in (\ref{section2:measure}) for $1\leq k \leq N$, that is we use the joint measure density function $m_k$ to replace $m$ and use $\textbf{M}_k$ to replace $\textbf{M}$. $\textbf{ME}_0$ is the finite measure expectation when no abrupt change occurs, which is the 1-order moment based on the pre-change joint measure density function $m_{0}$. That is
\begin{equation}
\textbf{ME}_{0}(\xi) = \int \xi d\textbf{M}_0
\end{equation}
where $\xi$ is a random variable, $\textbf{M}_0$ is the pre-change joint measure. Furthermore, if $\xi \in \mathfrak{F}_{n}=\sigma\{X_{k}, 0\leq k\leq n\}$ for $1\leq n\leq N$, then
$$
\textbf{ME}_{k}(\xi)=\begin{cases}
\textbf{ME}_{0}\Big(\xi\prod_{j=k}^{n}\Lambda_j\Big)& \text{$n \geq k$}\\
\textbf{ME}_{0}(\xi) & \text{$n < k$}
\end{cases}, \,\,\,\,\,\,\,\,\,    1 \leq k \leq N
$$
where $\Lambda_j=\frac{m_{1j}(X_j|X_{j-1}, ... , X_0)}{m_{0j}(X_j|X_{j-1}, ... , X_0)}$.

\begin{remark}
	It is easy to know that the measure expectation $\textbf{ME}_0$ and the expectation $\textbf{E}_0$ have similar properties.
\end{remark}

Let $T \in \mathfrak{T}_N$ be a sequential test defined as follows
$$
\begin{aligned}
T &= \min\{ 1\leq k\leq N+1: \, Y_k\geq l_k(c)\}\\
Y_k&=(Y_{k-1}+w_k)\Lambda_k, \Lambda_k=\frac{m_{1k}(X_k|X_{k-1}, ... , X_0)}{m_{0k}(X_k|X_{k-1}, ... , X_0)}
\end{aligned}
$$
for $1 \leq k \leq N$, and $Y_0:=0, Y_{N+1}:=Y_N$. $\mathfrak{T}_N$ is a set of all of the sequential tests satisfying $1\leq T\leq N+1$ and $ \{T\leq n\} \in  \mathfrak{F}_{n}=\sigma\{X_{k}, 0\leq k\leq n\}$ for $1\leq n\leq N$. Following the idea of the randomization probability of the change time and the definition of describing the average detection delay proposed by Moustakides (2008), we define a performance measure $\mathcal{J}_{N}(\cdot)$ which is based on the measure expectations \{$\textbf{ME}_k, 0 \leq k \leq N$\} to evaluate the detection performance of each sequential test $T\in \mathfrak{T}_N$ as follows
\begin{equation}
\label{section2.1:performance}
\mathcal{J}_{N}(T)=\frac{\sum_{k=1}^{N+1}\textbf{ME}_{k}(w_k(T-k)^+)}{\sum_{k=1}^{N+1}\textbf{ME}_{0}(v_kI(T \geq k))}=\frac{\sum_{k=1}^{N}\textbf{ME}_{k}(w_k(T-k)^+)}{\textbf{ME}_{0}(\sum_{k=1}^T v_k)}
\end{equation}
where $w_k(Y_{k-1}), v_k(Y_{k-1})$ are weights, which are functions of $Y_{k-1}$.

The widely used CUSUM test statistics can be written as
\begin{equation}
Y_k=\max \{1, Y_{k-1}\}\Lambda_k, \Lambda_k=\frac{m_{1k}(X_k|X_{k-1}, ... , X_0)}{m_{0k}(X_k|X_{k-1}, ... , X_0)}
\end{equation}
for $1\leq k \leq N$, and $Y_0:=0, Y_{N+1}:=Y_{N}$. This is a special case of $w_k=(1-Y_{k-1})^+$.

Note that $Y_k=\Lambda_k$ when $Y_{k-1} < 1$, or, $(1-Y_{k-1})^+>0$, that is, we will newly detect the change of networks starting from $k$. So it is reasonable to let the weight of the detection delay $w_k=(1-Y_{k-1})^+$. Actually, the numerator and denominator of the performance measure $\mathcal{J}_{N}(T)$ in (\ref{section2.1:performance}) can be regarded as the generalized out-of-control average run length $(ARL_1)$ and in-control average run length $(ARL_0)$. For the numerical simulations in Section \ref{section3}, we choose the CUSUM test statistics with $w_k=(1-Y_{k-1})^+, v_k=1$.  

\subsection{Optimal sequential tests} \label{section2.2}
In order to give the optimal network sequential test, we first give a definition about the optimal criterion of any sequential tests for finite network observations under the finite measure space. After defining a new performance measure under the finite measure space, we give the optimal sequential test, which is related to the control limits \{$l_k(c), 0 \leq k \leq N+1$\}. In this section, we discuss the optimal control limits \{$l_k(c), 0 \leq k \leq N+1$\} under the performance measure we defined.

\begin{definition}
	\label{definition1}
	A sequential test $T^* \in \mathfrak{T}_N$ with $\textbf{ME}_{0}(\sum_{k=1}^{T^*}v_k) \geq \gamma$ is  optimal under the performance measure $\mathcal{J}_{N}(T)$ if
	\begin{equation}
	\inf_{T\in \mathfrak{T}_N,\,\,\,\textbf{ME}_{0}(\sum_{k=1}^{T}v_k)\geq \gamma }\mathcal{J}_{N}(T)=\mathcal{J}_{N}(T^*)
	\end{equation}
	where $\gamma$ satisfies $\textbf{ME}_{0}(v_1)<\gamma <\textbf{ME}_{0}(\sum_{k=1}^{N+1}v_k)$. And note that $\textbf{ME}_{0}(\sum_{k=1}^{T}v_k)\geq \gamma$ means that the test should not exceed a certain false alarm rate, which is $\gamma$ here.
\end{definition}

Before giving the optimal sequential test, we need to give the definition of the corresponding conditional measure expectation w.r.t. (with respect to) $\textbf{M}_0$, which is also similar to the conditional expectation in the probability space. More details are as follows:

Let $\xi$ be an integrable random variable on $(\Omega, \mathcal{F}, \textbf{M})$. Let $\mathcal{A}$ be a sub-$\sigma$-field of $\mathcal{F}$. The conditional measure expectation of $\xi$ given $\mathcal{A}$ under the finite measure $\textbf{M}_0$, denoted by $\textbf{ME}_0(\xi|\mathcal{A})$, is the a.s.-unique random variable satisfying the following two conditions:

(i) $\textbf{ME}_0(\xi|\mathcal{A})$ is measurable in $(\Omega, \mathcal{A})$;

(ii) $\int_A \textbf{ME}_0(\xi|\mathcal{A})d\textbf{M}_0 = \int_A \xi d\textbf{M}_0$ for any $A \in \mathcal{A}$.

\begin{remark}
	It is easy to obtain that the conditional measure expectation $\textbf{ME}_0(\cdot|\mathcal{A})$ and the conditional expectation $\textbf{E}_0(\cdot|\mathcal{A})$ have similar properties. In the proofs of theorems, we need to use these properties.
\end{remark}

Motivated by Chow, Robbins, and Siegmund (1971), we present the nonnegative random  dynamic control limits $\{l_k(c), 0\leq k\leq N+1\}$ defined by the following recursive equations
\begin{equation}
\label{section2.2:dynamic}
\begin{aligned}
l_{N+1}(c)&=0,\,\,\, l_N(c)=cv_{N+1}  \\
l_{k}(c)&=cv_{k+1}+\textbf{ME}_{0}\Big([l_{k+1}(c)-Y_{k+1}]^+|\mathfrak{F}_{k}\Big)
\end{aligned}
\end{equation}
for $0\leq k\leq N-1$, where $c>0$ is a constant.

By using the test statistics $Y_n, 1\leq n\leq N+1,$ and the  control limits $l_k(c), 1\leq k\leq N+1$, we define a sequential test $T^*(c, N)$  as follows
\begin{equation}
T^*(c, N)=\min\{ 1\leq k\leq N+1: \, Y_k\geq l_k(c)\}.
\end{equation}
It is clear that $T^*(c, N)\in \mathfrak{T}_N$,  $l_k(c)\geq cv_{k+1}$ and $ l_k(c) \in \mathfrak{F}_{k}$ for $0\leq k\leq N$. The positive number $c$ can be regarded as an adjustment coefficient of the dynamic control limits, since $l_k(c), 0\leq k\leq N-1$ are increasing on $c\geq 0$ with $l_k(0)=0$ and $\lim_{c\to \infty}l_k(c)=\infty$ for $v_{k+1}>0$.

The following theorem shows that the sequential test $T^*(c, N)$ is optimal under the performance measure $\mathcal{J}_{N}(T)$.

\begin{theorem}
	\label{theorem1}
	Let $\gamma$ be a positive number satisfying $\textbf{ME}_0(v_1)<\gamma <\sum_{k=1}^{N+1}\textbf{ME}_0(v_k)$. There exists a positive number $c_{\gamma}$
	such that $T^*(c_{\gamma}, N)$ is optimal in the sense of (\ref{section2.1:performance}) with $\textbf{ME}_0(\sum_{k=1}^{T^*(c_{\gamma}, N)}v_k)=\gamma$; that is,
	\begin{equation}
	\inf_{T\in \mathfrak{T}_N, \,\, \textbf{ME}_{0}(\sum_{k=1}^{T}v_k)\geq \gamma }\mathcal{J}_{N}(T)=\mathcal{J}_{N}(T^*(c_{\gamma}, N)).
	\end{equation}
	In particular,  if $T\in \mathfrak{T}_N$ satisfies  $T\neq  T^*(c_{\gamma}, N)$, that is, $\textbf{M}_{0}( T\neq  T^*(c_{\gamma}, N))>0$ and  $\textbf{ME}_{0}(\sum_{k=1}^{T}v_k)=\gamma$, then
	\begin{equation}
	\mathcal{J}_{N}(T)>\mathcal{J}_{N}(T^*(c_{\gamma}, N)).
	\end{equation}
\end{theorem}
	
The proof of Theorem \ref{theorem1} is in the Appendix. Here, the random dynamic control limits $\{l_k(c), 0\leq k\leq N+1\}$ of the optimal control chart $T^*(c, N)$ can be called the optimal dynamic control limits. From Theorem \ref{theorem1}, we can construct an optimal control chart $T^*(c, N)$ for each of the performance measure $\mathcal{J}_{N}(T)$ corresponding to the weights $w_k, v_k$.

Next, we consider the Markov process. Since $\textbf{ME}_{0}([l_{k+1}(c)-Y_{k+1}]^+|\mathfrak{F}_{k})$ and $v_{k+1}$ are measurable with respect to $\mathfrak{F}_{k}$, there are functions $h_k=h_k(c, x_0, x_1, ..., x_k)$, $0\leq k\leq N-1$  and $v_{k}=v_{k}(x_0, x_1, ...,x_{k-1})$, $0\leq k\leq N-1$ such that
$$
h_k=h_k(c, x_0, x_1, ..., x_k)=\textbf{ME}_{0}([l_{k+1}(c)-Y_{k+1}]^+|X_k=x_k, X_{k-1}=x_{k-1}, ..., X_0=x_0)
$$
for $0\leq k\leq N-1$. Therefore, the optimal control limits \{$l_k(c), 0\leq k\leq N+1\}$ in (\ref{section2.2:dynamic}) can be written as
$$
l_k(c)=cv_{k+1}(X_0, X_1, ...,X_k) +h_k(c, X_0, X_1, ..., X_k)
$$
for $0\leq k\leq N$.

Note that the optimal control limits \{$l_k(c), 0\leq k\leq N+1\}$ depend on the sequence of networks, $ X_0, X_1, ..., X_k$ for $0\leq k\leq N$. Let the sequence of networks $\{X_k, \, 0\leq k\leq N\}$ be at most a $p$-order Markov process; that is, both the pre-change network observations $X_0, X_1, ..., X_{k-1}$ and the post-change network observations $X_{k}, ..., X_{N}$ are $i$-order and $j$-order Markov processes with transition measure density functions $m_{\textbf{0}n}(x_{n}|x_{n-1},..., x_{n-i})$ and  $m_{\textbf{1}s}(x_{s}|x_{s-1},..., x_{s-j})$, respectively, where $p=\max\{i,\,j\}$ and
$$
\begin{aligned}
m_{\textbf{0}n}(x_{n}|x_{n-1},..., x_{n-i}) & = m_{\textbf{0}n}(x_{n}|x_{n-1},..., x_{n-i},...,x_0)\\
m_{\textbf{1}s}(x_{s}|x_{s-1},..., x_{s-j}) & = m_{\textbf{1}s}(x_{s}|x_{s-1},..., x_{s-j},...,x_0)
\end{aligned}
$$
for $i\leq n$ and $j\leq s$. We know that $p=0$ means the sequence of networks are mutually independent. The definition of Markov process in the finite measure space is also similar.

\begin{theorem}
	\label{theorem2}
	Let the sequence of networks $\{X_k, \, 0\leq k\leq N\}$ be at most a $\textbf{p}$-order Markov process for $p\leq N$ and the weights $w_{k}=w_{k}(Y_{k-1})$
	and $v_{k}=v_{k}(Y_{k-1})$ for $1\leq k\leq N$. Then
	
	$(\romannumeral1)$ Let $p\geq 1$. The optimal control limits $\{l_k(c), \, 0\leq k\leq N\}$ can be written as
	$$
	l_k(c)=cv_{k+1}(Y_k)+\textbf{ME}_{0}\Big([l_{k+1}(c)-(Y_k+w_{k+1}(Y_k))\Lambda_{k+1}]^+|Y_k, X_k, X_{k-1},...,X_{k-p+1}\Big)
	$$
	for $0\leq k\leq N$. 
	
	$(\romannumeral2)$ Let  $p=0$ and $l_k(c)=l_k(c, Y_k)$, then we have
	$$
	l_k(c)=l_k(c, Y_k)=cv_{k+1}(Y_k)+\textbf{ME}_{0}\Big([l_{k+1}(c, Y_{k+1})-(Y_k+w_{k+1}(Y_k))\Lambda_{k+1}]^+|Y_k\Big)
	$$
	for $0\leq k\leq N$. 
\end{theorem}

Theorem \ref{theorem2} is a special case of Theorem \ref{theorem1} when the network observation sequence is a Markov process and the at most order of this Markov process $p$ is fixed.

\section{Comparison and analysis of simulation results} \label{section3}
Suppose that there is a finite sequence of networks and $N=60$. In Section \ref{section3.1}, we give an example of independent network observations by comparing four different sequential tests, $T_{cons}$ with constant control limits, $T_{de}$ with linear decrease control limits, $T_{in}$ with linear increase control limits, and $T^*(c)$ with optimal control limits $l_k(c)$ which are defined in this paper. Then we give an example of dependent network observations, a 1-order Markov process and also compare the four different sequential tests in Section \ref{section3.2}. All of the above are based on CUSUM-type statistics and the number of nodes $d=10$.

Note that the traditional probability space is a special case, so we can use the probability space when the normalization coefficient is computable for a network distribution. Hence we also give examples in Section \ref{section3.3} and Section \ref{section3.4} of independent and dependent Erd\H{o}s-R\'{e}nyi network sequences, respectively. Meanwhile, we can compare the similarities and differences between the finite measure space and the probability space through the results of the simulation.

In simulation, the performance measure $J_N(T)$ with $w_k=(1-Y_{k-1})^+$ and $v_k=1$ for $1 \leq k \leq N$ is a criterion. The smaller the value of $J_N(T)$, the better the performance of a sequential test. The ratio of $\textbf{ME}_0(T)$ to $\textbf{E}_0(T)$ is a deterministic and unknown constant for a sequential test $T$. $\textbf{ME}_0(T)$ or $\textbf{E}_0(T)$ is the false alarm rate. When comparing the performance of different sequential tests, we set a common false alarm rate, i.e.  $\textbf{ME}_0(T)$ or $\textbf{E}_0(T)$ and compare the value of $J_N(\cdot)$. In the third column of tables, we replace $\textbf{ME}_0(T)$ with $\textbf{E}_0(T)$. Note that the four tests all use adjustment coefficient $c$ to guarantee $\textbf{E}_0(T) \approx 40$, which is denoted by $c_\gamma$ in tables. 

For a finite network observation sequence, it is obvious that $1 \leq \textbf{E}_0(T) \leq N$, thus we can choose any value of $\textbf{E}_0(T) \in (1,N)$ in theory. Here we choose $\textbf{E}_0(T) \approx 40$, a suitable value to make the false alarm rate moderate. With the definition of a sequential test 
$$T(c, N)=\min\{ 1\leq k\leq N+1: \, Y_k\geq l_k(c,Y_k)\},$$
we can conclude that $c$ controls the values of the dynamic control limits, i.e. each value of $c$ corresponds to a value of $ARL_0(c)=\textbf{E}_0(T)$. That is, for a pre-determined false alarm rate, the adjustment coefficients $c$ for the four tests can be determined uniquely. Actually, for the four different sequential tests, the control limits are all increase functions of $c$, so $\textbf{ME}_0(T)$, i.e. $\textbf{E}_0(T)$, is also continuous and increasing on $c$, thus the values of $c_\gamma$'s for the four tests are unique.

Next we describe the simulation ideas. For any test $T$, $\textbf{E}_0(T)$ can be estimated by the average of repeated experiments. In each repetition, we first draw samples from the known pre-change measure density function. Once the samples are collected, we can get the statistics $Y_k, 1\leq k \leq N$ and compute the optimal dynamic control limits with the recursive equations for Markov process as follows
$$
\begin{aligned}
l_{N+1}(c)&=0,\,\,\, l_N(c)=cv_{N+1}  \\
l_{k}(c, Y_k)&=cv_{k+1}+\textbf{ME}_{0}\Big([l_{k+1}(c,Y_{k+1})-Y_{k+1}]^+|Y_{k}\Big)
\end{aligned}
$$
for $0\leq k\leq N-1$, where $c>0$ is a constant. Note that, the key point of computing the optimal dynamic control limits is the estimation of the functions $\textbf{ME}_{0}(\cdot{|Y_k})$, i.e. $l_{k}(c,Y_k)$ for $0\leq k\leq N-1$. Through the definition of a sequential test 
$$
T(c, N)=\min\{ 1\leq k\leq N+1: \, Y_k\geq l_k(c, Y_k)\},
$$
we can obtain the value of $T$. Finally, by repeating the above process, we can get the estimation of $\textbf{E}_0(T)$. The explanation of the simulation results is at the end of this section.

\subsection{Comparison of simulation of independent network observation sequences in the finite measure space} \label{section3.1}
Let $\{X_k, 0\leq k\leq 60\}$ be an independent network observation sequence with a pre-change measure density function $m_0(X_k) = \exp(-2 \times edges(X_k))$ and a post-change measure density function $m_1(X_k) = \exp(-2.2 \times edges(X_k))$ for $0\leq k\leq 60$. The function $edges$ is a term used in Exponential Random Graphical Model, which means the number of edges in a network, that is
$$
edges(X)=\begin{cases}
\frac{1}{2}\sum_{i,j=1}^d X_{ij}& \text{X is undirected}\\
\sum_{i,j=1}^d X_{ij}& \text{X is directed}
\end{cases}.
$$
There are many terms in ERGM and many linear combination of these terms. The example we give is the simplest. The ratio $\Lambda_k$ of the pre-change and post-change measure density functions $m_0(x)$ and $m_1(x)$ is
$$
\Lambda_k=\frac{m_1(X_k)}{m_0(X_k)}=e^{-0.2 \times edges(X_k)}
$$
for $1\leq k\leq 60$. We compare the performance of four tests, the first one is $T_{cons}$ with constant control limits
$$T_{cons}=\min\{ 1\leq k\leq N+1: \, Y_k\geq c\}.$$
The second test, $T_{de}$ denotes a test with decrease control lmits as follows
$$
\begin{aligned}
T_{de}&=\min\{ 1\leq k\leq N+1: \, Y_k\geq C(k)\},\\
C(k)&=c+0.005*(61-k).
\end{aligned}
$$
$T_{in}$ denotes a test with increase control lmits as follows
$$
\begin{aligned}
T_{in}&=\min\{ 1\leq k\leq N+1: \, Y_k\geq C(k)\},\\
C(k)&=c+0.005*(k+1).
\end{aligned}
$$
The last one is the optimal test $T^*(c)$ with $T^*(c)=\min\{ 1\leq k\leq N+1: \, Y_k\geq l_k(c)\}$. 

The simulation result for an independent and finite ERGM network observation sequence is shown in Table \ref{table1}. The result is based on $10^5$ repititions.

\begin{center}
	\label{table1}
	\textbf{Table 3.1}~~Simulation of four tests for an independent and finite ERGM network observation sequence.\\
	\setlength{\tabcolsep}{10mm}{
	\begin{tabular}{c|ccc}
		\hline\hline
		Sequential tests& $c_\gamma$ & $\textbf{E}_0$ & $J_N(T)$ \\
		\hline\hline
		$T_{cons}$ & 0.3679 & 41.8339 & 0.131814\\
		$T_{de}$ & 0.2165 & 40.0626 & 0.130858\\
		$T_{in}$ & 0.3277 & 39.7875 & 0.132746\\
		$T^*(c)$ & 0.1422 & 40.1929 & 0.130419\\
		\hline\hline
		\end{tabular}}
		\end{center}

\subsection{Comparison of simulation of dependent network observation sequences in the finite measure space} \label{section3.2}
Let $\{X_k, 0\leq k\leq 60\}$ be a dependent network observation sequence. And the measure density function of $X_0$ is $m(X_0)=1$, that is $X_0$ can been as a uniform distributed network. We note that if there is no cross term of $X_k$ and $X_{k-1}$ in $m(X_k|X_{k-1})$, then the network observations are independent obviously. So, here we let the pre-change conditional measure density function $m_0(X_k|X_{k-1}) = \exp(-0.08 \times edges(X_{k-1}) \times edges(X_k))$ and the post-change conditional measure density function $m_1(X_k|X_{k-1}) = \exp(-0.10 \times edges(X_{k-1}) \times edges(X_k))$ for $1\leq k\leq 60$. Then the ratio $\Lambda_k$ of the pre-change and post-change conditional measure density functions is
$$
\Lambda_k=\frac{m_1(X_k|X_{k-1})}{m_0(X_k|X_{k-1})}=e^{-0.02*edges(X_{k-1})*edges(X_k)}
$$
for $1\leq k\leq 60$. We compare the performance of four tests which are the same as those in Section \ref{section3.1}.

The simulation result of a dependent and finite ERGM network observation sequence is shown in Table \ref{table2}. The result is based on $10^4$ repititions. In this section, we only choose 1-order Markov process of network data.

\begin{center}
	\label{table2}
	\textbf{Table 3.2}~~Simulation of four tests for a dependent and finite ERGM network observation sequence.\\
	\setlength{\tabcolsep}{10mm}{
		\begin{tabular}{c|ccc}
			\hline\hline
			Sequential tests& $c_\gamma$ & $\textbf{E}_0$ & $J_N(T)$ \\
			\hline\hline
			$T_{cons}$ & 0.725 & 40.0135 & 0.077988\\
			$T_{de}$ & 0.700 & 39.8105 & 0.079665\\
			$T_{in}$ & 0.735 & 39.9017 & 0.080043\\
			$T^*(c)$ & 0.1295 & 39.9192 & 0.073666\\
			\hline\hline
		\end{tabular}}
	\end{center}
	
\subsection{Comparison of simulation of independent network observation sequences in the probability space} \label{section3.3}
The purpose of this section and next section is to compare the difference between the finite measure space and the probability space. Let $\{X_k, 0\leq k\leq 60\}$ be an independent network observation sequence. Before the change-point, the networks are Erd\H{o}s-R\'{e}nyi networks with link probability 0.5, and the link probability equals to 0.6 after the change-point. We compare the performance of four tests, the first one is $T_{cons}$ with constant control limits
$$T_{cons}=\min\{ 1\leq k\leq N+1: \, Y_k\geq c\}.$$
The second test, $T_{de}$ denotes a test with decrease control lmits as follows
$$
\begin{aligned}
T_{de}&=\min\{ 1\leq k\leq N+1: \, Y_k\geq C(k)\},\\
C(k)&=c+0.05*(61-k).
\end{aligned}
$$
$T_{in}$ denotes a test with increase control lmits as follows
$$
\begin{aligned}
T_{in}&=\min\{ 1\leq k\leq N+1: \, Y_k\geq C(k)\},\\
C(k)&=c+0.05*(k+1).
\end{aligned}
$$
The last one is the optimal test $T^*(c)$ with $T^*(c)=\min\{ 1\leq k\leq N+1: \, Y_k\geq l_k(c)\}$.

The simulation result of an independent and finite Erd\H{o}s-R\'{e}nyi network observation sequence is shown in Table \ref{table3}. The result is based on $10^5$ repititions.

\begin{center}
	\label{table3}
	\textbf{Table 3.3}~~Simulation of four tests for an independent and finite Erd\H{o}s-R\'{e}nyi network observation sequence.\\
	\setlength{\tabcolsep}{10mm}{
		\begin{tabular}{c|ccc}
			\hline\hline
			Sequential tests & $c_\gamma$ & $\textbf{E}_0$ & $J_N(T)$ \\
			\hline\hline
			$T_{cons}$ & 13.780 & 41.3517 & 1.240001 \\
			$T_{de}$ & 11.750 & 40.3647 & 1.235732  \\
			$T_{in}$ & 12.875 & 40.3902 & 1.267431  \\
			$T^*(c)$ & 2.050 & 39.9940 & 1.230929\\
			\hline\hline
		\end{tabular}}
	\end{center}

\subsection{Comparison of simulation of dependent network observation sequences in the probability space} \label{section3.4}
Let $\{X_k, 0\leq k\leq 60\}$ be a dependent network observation sequence. $X_0$ is a Erd\H{o}s-R\'{e}nyi network with the link probability 0.2. And given $X_{k-1}$, $X_{k}$ is a Erd\H{o}s-R\'{e}nyi network. Before the change-point, the link probability is $\frac{edges(X_{k-1})}{d(d-1)/2}$, and the link probability equals to $\frac{edges(X_{k-1})}{d(d-1)}$ after the change-point. 

We also compare the performance of four tests which are the same as those in Section \ref{section3.3}. The simulation result of a dependent and finite Erd\H{o}s-R\'{e}nyi network observation sequence is shown in Table \ref{table4}. The result is based on $10^4$ repititions.
\begin{center}
	\label{table4}
	\textbf{Table 3.4}~~Simulation of four tests for a dependent and finite Erd\H{o}s-R\'{e}nyi network observation sequence.\\
	\setlength{\tabcolsep}{10mm}{
		\begin{tabular}{c|ccc}
			\hline\hline
			Sequential tests& $c_\gamma$ & $\textbf{E}_0$ & $J_N(T)$ \\
			\hline\hline
			$T_{cons}$ & 8.850 & 39.9266 & 1.75347\\
			$T_{de}$ & 7.000 & 39.9224 & 1.676088\\
			$T_{in}$ & 8.410 & 40.2549 & 1.778835\\
			$T^*(c)$ & 4.000 & 40.0501 & 1.293891\\
			\hline\hline
		\end{tabular}}
	\end{center}
	
Finally, we give a explanation about the numerical simulation results in the above four tables. The first two tables are for general independent and dependent ERGM. The last two tables are for specific Erd\H{o}s-R\'{e}nyi random networks. The first column shows that we compare four representative different sequential tests, the adjustment coefficients $c_\gamma$'s are shown in the second column to guarantee that $\textbf{ME}_0(\cdot) \approx \gamma$, we use $\textbf{E}_0(\cdot)$ to replace it. With the commonly pre-determined false alarm rate, we compare the performance measure $J_N(T)$. Through the Definition \ref{definition1}, we can conclude that the smaller the value of $J_N(T)$, the better the performance of a sequential test. In this view, the fourth column is the focus of our attention. We find that the optimal sequential tests in four tables all have the smallest $J_N(T)$, which is consistent with Theorem \ref{theorem1}. Note that the method in this paper can be applied in many networks, as long as we know the distribution of the networks either in the probability space or in a general finite measure space. The sequence of networks in Section \ref{section3.2} and Section \ref{section3.4} can be $p$-order, while in simulation we choose 1-order networks. As we try to make sure $\textbf{E}_0(\cdot) \approx 40$, but in Table \ref{table1} and Table \ref{table3}, for the constant control limits, $\textbf{E}_0(\cdot)=41.8339, 41.3517$, that is because if we adjust the constant $c$ a little bit, the corresponding $\textbf{E}_0(\cdot)$ is far from 40. 

There are many different control limits, in this paper we choose relatively representative control limits. The usual CUSUM control chart we use in quality detection is constant control limits, denoted by $T_{cons}$ in this paper. From the simulation, we can conclude that under the performance measure $J_N(T)$, constant control limits maybe not good enough. Somtimes the increase control limits or the decrease control limits have better performance, which motivates efforts to find the optimal dynamic control limits. The optimal sequential test we defined is the guidance for how to give the suitable or optimal dynamic control limits. Indeed, the optimal control limits in simulation are a.s. decrease. Here, in above tables, we compute the value of $J_N(\cdot)$ for different sequential tests when they have a common generalized $ARL_0=\textbf{ME}_0(T)$ (the denominator of $J_N(T)$). The essence of the performance measure $J_N(\cdot)$ is to compare the value of generalized $ARL_1$, i.e. the numerator of $J_N(\cdot)$ when the value of generalized $ARL_0$, i.e. the denominator of $J_N(\cdot)$ is the same.

\section{Concluding remarks} \label{section4}
For the models such as ERGM and probabilistic graphical model, the normalization coefficient is difficult to compute for large scale networks. There are many methods to approximate the coefficient, which needs complex calculations and the accuracy is not high. We construct a finite measure space and define corresponding measure ratio statistics for change detection problem, which is a method to overcome the difficulty. However, it does not mean that the likelihood ratio statistics do not work. Measure ratio statistics can also show the difference between pre-change and post-change distribution. We can also give the optimal sequential test under the finite measure space.

When the pre-change and post-change measure density functions are known, the change detection problem can be solved. However, if the post-change measure density function is unknown, the problem is more difficult and worthy of studying as a new research topic. If the number of parameters in ERGM is very large, it may suffer dimension curse, which is also worthy to be considered. Dimension curse also exists in the problem of change-point monitoring, which can mislead the judgement (Donoho and Johnstone 1994). Even if there is no change, the superposition of high-dimensional noise errors will produce a high false alarm rate.
\end{spacing}

\newpage
\noindent \textbf{APPENDIX : PROOFS OF THEOREMS}

By using similar method given in Han, Tsung, and Xian (2019), we can prove Theorem \ref{theorem1} and Theorem \ref{theorem2}. The expectation and conditional expectation will be replaced by measure expectation and conditional measure expectation in the following.
\setcounter{equation}{0}
\renewcommand\theequation{A. \arabic{equation}}

{\bf Proof of Theorem \ref{theorem1}.} Han, Tsung, and Xian (2019) considered the probability space and probability expectations, we extend to a finite measure space with measure expectations. Similarly, $(T-k)^+=\sum_{m=k+1}^{N+1}I(T \geq m)$ for $T\in \mathfrak{T}_N$, $I(T \geq m) \in \mathfrak{F}_{m-1}$, we have
\begin{equation}
\label{appendix:transform}
\begin{aligned}
\sum_{k=1}^{N}\textbf{ME}_{k}\Big(w_k(T-k)^+\Big)&=\textbf{ME}_{0}\Big(\sum_{k=1}^{N}\sum_{m=k+1}^{N+1}w_kI(T\geq m)\prod_{j=k}^{m-1}\Lambda_j\Big)\\
&=\textbf{ME}_{0}\Big(\sum_{m=1}^{N}Y_{m}I(T\geq m+1)\Big)=\textbf{ME}_{0}\Big(\sum_{m=1}^{T}Y_{m-1}\Big),
\end{aligned}
\end{equation}
for all $T\in \mathfrak{T}_N$.

Let $T^*=T^*(c)=T^*(c, N)$ and
\begin{equation}
\label{appendix:xi}
\xi_n=\sum_{k=1}^n(Y_{k-1}-cv_k),
\end{equation}
where $c>0$. This is divided into two steps to complete the proof of Theorem \ref{theorem1}.

\textbf{Step I.} Show that
\begin{equation}
\label{appendix:Step1}
\textbf{ME}_{0}(\xi_T) \geq \textbf{ME}_{0}(\xi_{T^*})
\end{equation}
for all $T\in \mathfrak{T}_N$ and the strict inequality of (\ref{appendix:Step1}) holds for all $T\in \mathfrak{T}_N$ with $T\neq T^*$.

To prove (\ref{appendix:Step1}), by Lemma 3.2 in Chow, Robbins and Siegmund (1971), we extend to measure space, that is we only need to prove the following two inequalities:
\begin{equation}
\label{appendix:step1-1}
\textbf{ME}_{0}(\xi_{T^*}|\mathfrak{F}_{n}) \leq \xi_{n}\,\,\,\, \text{ on } \,\, \{ T^*> n\}
\end{equation}
and
\begin{equation}
\label{appendix:step1-2}
\textbf{ME}_{0}(\xi_{T}|\mathfrak{F}_{n}) \geq \xi_{n}\,\,\,\, \text{ on } \,\, \{T^*=n, \, T > n\}
\end{equation}
for each $n\geq 1$.

Using similar method as in Han, Tsung, and Xian(2019), we can verify that for $1\leq n\leq N$,
\begin{equation}
\label{appendix:step1-1process}
I(T^*> n)\textbf{ME}_{0}(\xi_{T^*}-\xi_{n})|\mathfrak{F}_{n})=I(T^*> n)(Y_n-l_n(c))<0,
\end{equation}
and
\begin{equation}
\label{appendix:step1-2process}
I(T^*= n)I(T>n)\textbf{ME}_{0}(\xi_{T}-\xi_{n})|\mathfrak{F}_{n}) \geq I(T^*=n)I(T>n)[Y_n-l_n(c)]\geq 0.
\end{equation}
From (\ref{appendix:step1-1process}) and (\ref{appendix:step1-2process}) it follows that (\ref{appendix:step1-1}) and (\ref{appendix:step1-2}) hold for $1 \leq n \leq N$. By (\ref{appendix:step1-1}) and (\ref{appendix:step1-2}), we know that the inequality in (\ref{appendix:Step1}) holds for all $T\in \mathfrak{T}_N$. Furthermore, from (\ref{appendix:step1-1process}) and (\ref{appendix:step1-2process}) it follows that the strict inequality in (\ref{appendix:Step1}) holds for all $T\in \mathfrak{T}_N$ with $T\neq T^*$.

\textbf{Step II.} Show that the test we defined in this paper is optimal under the performance measure $\mathcal{J}_{N}$.  

Let $k^*=\max\{2\leq k\leq N+1: \,\,\textbf{ME}_0(v_k)>0\}$, as in Han, Tsung, and Xian(2019), it follows that there is a positive number $c_{\gamma}$ such that
\begin{equation}
\textbf{ME}_{0}(\sum_{k=1}^{T^*(c_{\gamma})}v_k)=\sum_{k=1}^{k^*}\textbf{ME}_{0}(v_kI(T^*(c_{\gamma})\geq k))=\gamma.
\end{equation}

Let
$$
\tilde{c}_{\gamma}= \mathcal{J}_{N}(T^*(c_{\gamma}))=c_{\gamma}-\frac{\textbf{ME}_{0}(l_0(c_{\gamma}))}{\textbf{ME}_{0}(\sum_{m=1}^{T^*(c_{\gamma})}v_m)}.
$$
If $\mathcal{J}_{N}(T)\geq c_{\gamma}$, then $\mathcal{J}_{N}(T)\geq \tilde{c}_{\gamma}= \mathcal{J}_{N}(T^*(c_{\gamma}))$. If $\mathcal{J}_{N}(T)<  c_{\gamma}$, then,
by (\ref{appendix:transform}), (\ref{appendix:Step1}), and $\textbf{ME}_{0}(\sum_{m=1}^{T}v_m)\geq \gamma$,  we have
$$
\begin{aligned}
(\mathcal{J}_{N}(T)-c_{\gamma}) \gamma & \geq  [\frac{ \textbf{ME}_{0}(
\sum_{m=1}^{T}Y_{m-1})}{\textbf{ME}_{0}(\sum_{m=1}^{T}v_m)}-c_{\gamma}]\textbf{ME}_{0}(\sum_{m=1}^{T}v_m)\\
&=[\textbf{ME}_{0}(
\sum_{m=1}^{T}Y_{m-1})-c_{\gamma}\textbf{ME}_{0}(\sum_{m=1}^{T}v_m)]\\
&\geq  [\textbf{ME}_{0}( \sum_{m=1}^{T^*(c_{\gamma})}Y_{m-1})-c_{\gamma}\textbf{ME}_{0}(\sum_{m=1}^{T^*(c_{\gamma})}v_m)]\\
&=[\frac{ \textbf{ME}_{0}( \sum_{m=1}^{T^*(c_{\gamma})}Y_{m-1})}{\textbf{ME}_{0}(\sum_{m=1}^{T^*(c_{\gamma})}v_m)}-c_{\gamma}]\textbf{ME}_{0}(\sum_{m=1}^{T^*(c_{\gamma})}v_m)\\
&=[\mathcal{J}_{N}(T^*(c_{\gamma}))-c_{\gamma}]\gamma.
\end{aligned}
$$
This means that $\mathcal{J}_{N}(T)\geq \mathcal{J}_{N}(T^*(c_{\gamma}))$ for all $T\in \mathfrak{T}_N$ with $\textbf{ME}_{0}(\sum_{m=1}^{T}v_m)\geq \gamma$. The strict inequality of Theorem \ref{theorem1} comes from the strict inequality in (\ref{appendix:Step1}) when  $T\neq T^*(c_{\gamma})$ with $\textbf{ME}_{0}(\sum_{m=1}^{T}v_m)=\textbf{ME}_{0}(\sum_{m=1}^{T^*(c_{\gamma})}v_m)= \gamma$. This completes the proof of Theorem \ref{theorem1}.

{\bf Proof of Theorem \ref{theorem2}. } As $Y_k=(Y_{k-1}+w_{k}(Y_{k-1}))\Lambda_k$ for $1\leq k\leq N$, it follows that $(Y_k, X_k), 0\leq k\leq N,$ is at most a two-dimensional $p$-order Markov process in the finite measure space. Let $p\geq 1$. By the definition of the optimal control limits, we have
\begin{equation}
l_k(c)=cv_{k+1}(Y_k)+\textbf{ME}_{0}\Big([l_{k+1}(c)-(Y_k+w_{k+1}(Y_k))\Lambda_{k+1}]^+|Y_k, X_k,..., X_0\Big)
\end{equation}
for $0\leq k\leq  p-1$ and
\begin{equation}
\begin{aligned}
l_k(c)&=cv_{k+1}(Y_k)\\
&+\textbf{ME}_{0}\Big([l_{k+1}(c)-(Y_k+w_{k+1}(Y_k))\Lambda_{k+1}]^+|Y_k, X_k,..., X_{k-p+1}\Big)
\end{aligned}
\end{equation}
for $ p \leq k\leq N.$ Let $p=0$,  we have similarly
\begin{equation}
\begin{aligned}
l_k(c)&=l_k(c, Y_k)\\
&=cv_{k+1}(Y_k)+\textbf{ME}_{0}\Big([l_{k+1}(c, Y_{k+1})-(Y_k+w_{k+1}(Y_k))\Lambda_{k+1}]^+|Y_k\Big)
\end{aligned}
\end{equation}
for $0\leq k\leq N$.
\end{document}